\documentclass[a4paper]{article}	
\usepackage{hyperref}
\usepackage{graphicx}
\usepackage{amsfonts}
\usepackage{amssymb}
\usepackage{amsmath}
\usepackage{array}
\usepackage{authblk}
\usepackage{cite}
\usepackage{epstopdf}
\usepackage{float}
\newcommand{\PreserveBackslash}[1]{\let\temp=\\#1\let\\=\temp}
\newcolumntype{C}[1]{>{\PreserveBackslash\centering}p{#1}}
\newcolumntype{R}[1]{>{\PreserveBackslash\raggedleft}p{#1}}
\newcolumntype{L}[1]{>{\PreserveBackslash\raggedright}p{#1}}
\begin{document}
\setlength{\arraycolsep}{0pt}
\setlength{\baselineskip}{15pt}
\setlength{\parindent}{2em}
\title{\bf Some Sufficient Conditions on Pancyclic Graphs}
\author{Guidong Yu$^{1,2}$\thanks{Email: guidongy@163.com.
Supported by the Natural Science Foundation of China (No. 11871077), the NSF of Anhui Province (No.
1808085MA04), and the NSF of
Department of Education of Anhui Province (No. KJ2017A362).}, Tao Yu$^{1}$, Axiu Shu$^{1}$, Xiangwei Xia$^{1}$\\
  {\small \it  $1.$ School of Mathmatics and Computation Sciences, Anqing
Normal University, Anqing 246133,
China.}\\
{\small \it  $2.$ Basic Department, Hefei Preschool Education College, Hefei
230013, China.} }
\renewcommand*{\Affilfont}{\small\it}
\renewcommand\Authands{ and }
\maketitle
{\small
\noindent{\bf Abstract:}
A pancyclic graph is a graph that contains cycles of all possible lengths from three up to the number of vertices in the graph. In this paper, we establish some new sufficient conditions for a graph to be pancyclic  in terms of the edge number, the spectral radius and the signless Laplacian spectral radius of the graph.


\noindent {\bf MR Subject Classifications:} 05C50,15A18.

\section{Introduction}
In this paper, we use $G=(V(G),E(G))$ to denote a finite simple undirected graph with vertex set $V(G)=\{v_1,v_2,\cdots,v_n\}$ and edg set $E(G)$. Write by $m=|E(G)|$ the number of edges of the graph $G$. Let $v_i\in V(G)$, we denote by $d_i=d_{v_{i}}=d_G(v_i)$ the degree of $v_i$. Let $(d_1,d_2,\cdots,d_n)$ be the degree sequence of $G$, where $d_1\leq d_2\leq \cdots\leq d_n$. Denote by $\delta(G)$ or simply $\delta$ the minimum degree of $G$. The set of neighbours of a vertex $v$ in $G$ is denoted by $N_G(v)$. We use $G[X,Y]$ to denote a bipartite graph with bipartition $(X,Y)$. Let $K_n$ be a complete graph of order $n$ and $K_{m,n}$ be a complete bipartite graph with two parts having $m,n$ vertices, respectively. Let $G$ and $H$ be two disjoint graphs. The disjoint union of $G$ and $H$, denoted by $G+H$, is the graph with vertex set $V(G)\cup V(H)$ and edge set $E(G)\cup E(H)$. If $G_1=G_2=\cdots=G_k$, we denote $G_1+G_2+\cdots+G_k$ by $kG_1$. The join of $G$ and $H$, denoted by $G\vee H$, is the graph obtained from disjoint union of $G$ and $H$ by adding edges joining every vertex of $G$ to every vertex of $H$. Let $K_{n-1}+v$ denote the complete graph on $n-1$ vertices together with an isolated vertex $v$.

  The adjacency matrix of $G$ is defined to be a matrix $A(G)=[a_{ij}]$ of order $n$, where $a_{ij}=1$ if $v_i$ is adjacent to $v_j$, and $a_{ij}=0$ otherwise. The largest eigenvalue of $A(G)$, denote by $\mu (G)$, is called to be the spectral radius of $G$. Let $D(G)$ be the drgree diagonal matrix of $G$. The matrix $Q(G)=D(G)+A(G)$ is the signless Laplacian matrix of $G$. The largest eigenvalue of $Q(G)$, denoted by $q(G)$, is called to be the signless Laplacian spectral radius of $G$.

 A cycle (path) containing all vertices of a graph $G$ is called a Hamilton cycle (path) of $G$. A graph $G$ is
hamiltonian if it contains a Hamilton cycle. And $G$ is pancyclic if it contains cycles of every
length $l$, $3 \leq l \leq n$. Clearly, a bipartite graph is not pancyclic. A pancyclic graph is certainly
Hamiltonian, but the converse is not true. A cycle of length $l$ is called an $l-$cycle. The problem of deciding whether a graph is Hamiltonian is one of the most difficult classical problems in graph theory. Indeed, it is {\it NP-complete}.

Recently, the spectral theory of graphs has been applied to this problem. Firstly, Fiedler and Nikiforov \cite{21} gave tight conditions on spectral radius of a graph and its complement
for the existence of Hamiltonian paths and cycles. Next, Bo Zhou \cite{22} gave tight conditions on the
signless spectral radius of a graph complement for the existence of Hamiltonian paths and cycles. Yu and Fan \cite{3} established the spectral conditions for a graph to be
Hamilton-connected in terms of the spectral radius of the adjacency matrix or signless Laplacian matrix of
the graph or its complement. Lu, Liu and Tian \cite{23} gave sufficient conditions
for a bipartite graph to be Hamiltonian in terms of the spectral radius of the adjacency matrix of
the graph. Since then, many researchers have studied the
analogous problems under various spectral conditions; see \cite{4, 26, 24, 27, 9, 28, 5}. But there is no spectral sufficient conditions on pancyclic graphs.
In this paper, we first establish a new sufficient conditions for a graph to be pancyclic in terms of the edge number of the graph, then basing on edge number sufficient condition, we give (signless Laplacian) spectral radius sufficient conditions for a graph to be pancyclic.

\section{Preliminary}
We begin with some definitions. Given a graph $G$ of order $n$, a vector $X\in R^n$ is called to be defined on $G$, if there is a 1-1 map $\varphi$ from $V(G)$ to the entries of $X$; simply written $X_u=\varphi(u)$.

If $X$ is an eigenvector of $A(G)$ ($Q(G)$), then $X$ is defined naturally on $G$, i.e. $X_u$ is the entry of $X$ corresponding to the vertex $u$. One can find that
when $\lambda$ is an eigenvalue of $G$ corresponding to the eigenvector $X$ if and only if $X\neq0$,
$$\lambda X_v=\sum\limits_{u\in N_G(v)} X_u, \mbox{~for each vertex }v\in V(G).\eqno(2.1)$$
The equation (2.1) is called {\it eigen-equation} of $G$.
When $q$ is an signless Laplacian eigenvalue of $G$ corresponding to the eigenvector $X$ if and only if $X\neq0$, one can find that
$$[q-d_G(v)]X_v=\sum\limits_{u\in N_G(v)} X_u, \mbox{~for each vertex } v\in V(G).\eqno(2.2)$$
The equation (2.2) is called {\it signless Laplacian eigen-equation} of $G$.

\noindent{\bf Lemma 2.1}\cite{12} Let $G$ be a graph of order $n$ with degree sequence $d_1\le\cdots\le d_n$, if for all positive integers $k$ such that $d_k\le k\textless\frac{n}{2}$ and $d_{n-k}\ge{n-k}$, then $G$ is a pancyclic graph or bipartite graph.

\noindent{\bf Lemma 2.2}\cite{1} Let $G$ be a connected graph of order $n$ with $m$ deges. Then
\[\mu(G)\le\sqrt{2m-n+1,}\]
and the equality holds if and only if $G=K_n$ or $G=K_{1,n-1}$.

\noindent{\bf Lemma 2.3}\cite{3} Let $G$ be a graph of order $n$ with $m$ deges.
Then
\[q(G)\le\frac{2m}{n-1}+n-2.\]
If $G$ is connected, the equality holds if and only if $G=K_{1,n-1}$ or $G=K_n$. Otherwise, the equality holds if and only if $G$=$K_{n-1}+v$.

\section{Main Results}

{\bf Theorem 3.1}	Let $G$ be a connected graph on $n(\ge 5)$ vertices and $m$ deges with minimum degree $\delta(G)\ge2$. If
	\[m\ge\left(
		\begin{array}{c}
		n-2\\	
		2\\
		\end{array}
	\right)+4,\eqno(3.1)\]
then $G$ is a pancyclic graph unless $G$ is a bipartite graph or $G\in \mathbb{NP}_1$=$\{K_2\vee(K_{n-4}+2K_1), K_5\vee6K_1, K_3\vee(K_2+3K_1), K_3\vee(K_1+K_{1,4}), K_3\vee(K_2+K_{1,3}), (K_2\vee2K_1)\vee5K_1, K_4\vee5K_1, K_{1,2}\vee4K_1, K_2\vee(K_1+K_{1,3}), K_3\vee4K_1\}$.

\noindent{\bf Proof:}  Suppose that $G$ is neither a pancyclic graph nor a bipartite graph. By {\bf Lemma 2.1}, there exists an positive integer $k$ for $d_k\le k\textless\frac{n}{2}$, such that $d_{n-k}\le n-k-1$. Then we have
\begin{equation*}
\begin{aligned}
2m & =\sum\limits_{i=1}^k d_i+\sum\limits_{i=k+1}^{n-k} d_i+\sum\limits_{i=n-k+1}^n d_i\\
   & \le k^2+(n-2k)(n-k-1)+k(n-1)\\
   & =n^2-n+3k^2+(1-2n)k\\
   & =	
	2\left(
		\begin{array}{c}
		n-2\\	
		2\\
		\end{array}
	\right)+8-(k-2)(2n-3k-7),
\end{aligned}
\end{equation*}

\noindent thus
	\[m\le
	\left(
		\begin{array}{c}
		n-2\\	
		2\\
		\end{array}
	\right)+4-\frac{(k-2)(2n-3k-7)}{2}.\eqno(3.2)\]
Since
$\left(
		\begin{array}{c}
		n-2\\	
		2\\
		\end{array}
	\right)+4\le m \le
\left(
		\begin{array}{c}
		n-2\\	
		2\\
		\end{array}
	\right)+4-\frac{(k-2)(2n-3k-7)}{2}$, thus $(k-2)(2n-3k-7)\le0$. Next, we discuss in the follow two cases.

\noindent{\bf Case 1} Assume that $(k-2)(2n-3k-7)=0$, i.e., $k=2$ or $2n-3k-7=0$.

\noindent Then, $m=
\left(
		\begin{array}{c}
		n-2\\
		2\\
		\end{array}
	\right)+4$ and all inequalities in the above arguments should be equalities.

\noindent{\bf Case 1.1} If $k=2$, then $G$ is a graph with  $d_1=d_2=2$, $d_3=\cdots=d_{n-2}=n-3$, $d_{n-1}=d_n=n-1$. The two vertices of degree $n-1$ must be adjacent to every vertex, so they induce a $K_2$. The two vertices of degree $2$ are not adjacent to other vertices, so they induce a $2K_1$. For the remaining $n-4$ vertices of degree $n-3$. They must be adjacent to each other to make sure the requirement of the degree $n-3$, so they induce a $K_{n-4}$. By the above analysis, we can get the graph $G$ must be  $K_2\vee(K_{n-4}+2K_1)$.

\noindent{\bf Case 1.2} If $2n-3k-7=0$, then we can get $n\le13$ because $k\textless\frac{n}{2}$, and hence  $n=11$, $k=5$ ,or $n=8$, $k=3$. The corresponding permissible graphic sequence are  $(5,5,5,5,5,5,10,10,10,10,10)$, $(3,3,3,4,4,7,7,7)$, respectively.

For the degree sequence $(5,5,5,5,5,5,10,10,10,10,10)$. The five vertices of degree $10$ must be adjacent to every vertex, so they induce a $K_5$. The remaining six vertices now have degree $5$, so they induce a $6K_1$. Then the graph must be  $K_5\vee6K_1$. By the similar discussion, the degree sequence $(3,3,3,4,4,7,7,7)$ must be correspond to $K_3\vee(K_2+3K_1)$.

\noindent{\bf Case 2} Assume that $(k-2)(2n-3k-7)<0$, $i.e.,$ $k\ge3$ and $2n-3k-7\textless0$. In this case, we have $6\le2k<n\le13$.

\noindent{\bf Case 2.1} If $n=13$, then $k\le6$ and $2n-3k-7\textgreater0$.

\noindent{\bf Case 2.2} If $n=12$, then $k\le5$ and $2n-3k-7\textgreater0$.

\noindent{\bf Case 2.3} If $n=11$, then $k\le5$ and $2n-3k-7\ge0$.

\noindent{\bf Case 2.4} If $n=10$, then $k\le4$ and $2n-3k-7\textgreater0$.

\noindent{\bf Case 2.5} If $n=8$, then $k\le3$ and $2n-3k-7=0$.

The above five cases all contradict to $2n-3k-7\textless0$.

\noindent{\bf Case 2.6} If $n=9$, then $k\le4$. When $k=4$, then $d_4\le4$, $d_5\le4$ and we have $50\le\sum\limits_{i=1}^9 d_i\le52$ by (3.1) and (3.2). When $\sum\limits_{i=1}^9 d_i=50$, the degree sequence of $G$ is $(3,4,4,4,4,7,8,8,8)$ or $(4,4,4,4,4,6,8,8,8)$ or $(4,4,4,4,4,7,7,8,8)$, it is easy to see that $G=K_3\vee(K_1+K_{1,4})$ or $G=K_3\vee(K_2+K_{1,3})$ or $G=(K_2\vee2K_1)\vee5K_1$; When $\sum\limits_{i=1}^9 d_i=52$, the degree sequence of $G$ is $(4,4,4,4,4,8,8,8,8)$, it is easy to see that $G= K_4\vee5K_1$. When $k=3$, then $2n-3k-7=2\textgreater0$, contradiction with $2n-3k-7\textless0$.

\noindent{\bf Case 2.7} If $n=7$, then $k\le3$, and because $k\ge3$, i.e., $k=3$ . Thus $d_3\le3$, $d_4\le3$ and $28\le\sum\limits_{i=1}^7 d_i\le30$. When $\sum\limits_{i=1}^7 d_i=28$, the degree sequence of $G$ is $(3,3,3,3,5,5,6)$ or (3,3,3,3,4,6,6) or (2,3,3,3,5,6,6), it is easy to see that $G=K_{1,2}\vee4K_1$ or $G=K_2\vee(K_2+K_{1,2})$ or $G=K_2\vee(K_1+K_{1,3})$. When $\sum\limits_{i=1}^7 d_i=30$, the degree sequence of $G$ is $(3,3,3,3,6,6,6)$, it is easy to see that $G=K_3\vee4K_1$.

\begin{table}[H]
\caption{The maximum length cycle $(l(G))$ of $G$}
		\renewcommand\arraystretch{1.5}
		\centering
		\begin{tabular}{C{3cm}C{2cm}C{3cm}C{2cm}}\hline
			$G$&$l(G)$&$G$&$l(G)$\\\hline
$K_2\vee(K_{n-4}+2K_1)$&$C_{n-1}$&$K_3\vee4K_1$&$C_6$\\
			$K_{1,2}\vee4K_1$&$C_6$&$K_3\vee(K_2+3K_1)$&$C_7$\\
			$K_2\vee(K_1+K_{1,3})$&$C_6$&$K_3\vee(K_2+K_{1,3})$&$C_8$\\
			$K_2\vee(K_2+K_{1,2})$&$C_7$&$K_4\vee5K_1$&$C_8$\\
			$K_3\vee(K_1+K_{1,4})$&$C_7$&$K_5\vee6K_1$&$C_{10}$\\
                $(K_2\vee2K_1)\vee5K_1$&$C_8$\\\hline
		\end{tabular}
		\label{tab:1}
		\end{table}
In {\it Table 1}, $G=K_2\vee(K_2+K_{1,2})$ contains cycles of every length $l$, $3\le l \le7$, namely it is pancyclic graph, a contradiction. The other graphs in {\it Table 1} are neither pancyclic nor bipartite.

The proof is complete.$\blacksquare$

\noindent{\bf Theorem 3.2} Let $G$ be a connected graph on $n(\ge5)$ vertices with minimum degree $\delta(G)\ge2$. If
\[\mu(G)\ge\sqrt{n^2-6n+15},\]
then $G$ is a pancyclic graph unless $G$ is a bipartite graph.

\noindent{\bf Proof:} Suppose that $G$ with $m$ edges is neither a pancyclic graph nor a bipartite graph. Because $K_n$ is pancyclic and $\delta(K_{1,n-1})=1$. By {\bf Lemma 2.2},
\[\sqrt{n^2-6n+15}\le\mu(G)<\sqrt{2m-n+1},\]
then
\[m>\left(
		\begin{array}{c}
		n-2\\	
		2\\
		\end{array}
	\right)+4.\]

By {\bf Theorem 3.1}, we get $G\in \mathbb{NP}_1$=$\{K_2\vee(K_{n-4}+2K_1), K_5\vee6K_1, K_3\vee(K_2+3K_1), K_3\vee(K_1+K_{1,4}), K_3\vee(K_2+K_{1,3}), (K_2\vee2K_1)\vee5K_1, K_4\vee5K_1, K_{1,2}\vee4K_1, K_2\vee(K_1+K_{1,3}), K_3\vee4K_1\}$. According to calculation, when $G\in \{K_2\vee(K_{n-4}+2K_1), K_5\vee6K_1, K_3\vee(K_2+3K_1)\}$, \[m=\left(
		\begin{array}{c}
		n-2\\	
		2\\
		\end{array}
	\right)+4,\] a contradiction.
So, $G\in \mathbb{NP}_2$=$\{K_3\vee(K_1+K_{1,4}), K_3\vee(K_2+K_{1,3}), (K_2\vee2K_1)\vee5K_1, K_4\vee5K_1, K_{1,2}\vee4K_1, K_2\vee(K_1+K_{1,3}), K_3\vee4K_1\}$.

For $G=(K_2\vee2K_1)\vee5K_1$, let $X=(X_1,X_2,\cdots,X_9)^T$ be the eigenvector corresponding to $\mu (G)$, where $X_i(1\le i\le 5)$ correspond to the vertex of degree $4$, $X_i(6\le i\le 7)$ correspond to the vertex of degree $7$ and $X_i(8\le i\le 9)$ correspond to the vertex of degree $8$. Then by eigen-equation (2.1), we have
\begin{equation*}
\setlength{\arraycolsep}{5pt}
\left\{
  \begin{array}{l}
X_1=X_2=\cdots=X_5,X_6=X_7,X_8=X_9,\\
\mu(G)X_1=2X_6+2X_8,\\
\mu(G)X_6=5X_1+2X_8,\\
\mu(G)X_8=5X_1+2X_6+X_8.\\
  \end{array}
\right.
\end{equation*}
Transform above into a matrix equation $(A'(G)-\mu(G)I)X'=0$, where $X'=(X_1,X_6,X_8)^T$ and
\begin{equation*}
\setlength{\arraycolsep}{5pt}
A'(G)=\left[		
		\begin{array}{ccc}		
		0 & 2 & 2\\		
		5 & 0 & 2\\		
		5 & 2 & 1\\
		\end{array}
	\right].
\end{equation*}
Let $f(x) :=det(xI-A'(G))$, then $f(x)=x^3-x^2-24x-30$, and $\mu (G)$ is the largest root of $f(x)=0$. Through calculation, $\mu(G)=5.9150<\sqrt{9^2-
6\times9+15}$, a contradiction.
Using the same method, we get the spectral radius of the other graphs in $\mathbb{NP}_2$, showing in the following {\it Table 2}.
\begin{table}[H]
\caption{The spectral radius of $G$}
		\renewcommand\arraystretch{1.5}
		\centering
		\resizebox{\textwidth}{!}{
		\begin{tabular}{C{2.5cm}C{1cm}C{2.5cm}C{2.5cm}C{1cm}C{2.5cm}}\hline
			$G$&$\mu(G)$&$\sqrt{n^2-6n+15}$&  $G$&$\mu(G)$&$\sqrt{n^2-6n+15}$\\\hline
			$K_{1,2}\vee4K_1$&$4.2182$&$4.6904$&$K_3\vee(K_2+K_{1,3})$&$5.9612$&$6.4807$\\
			$K_2\vee(K_1+K_{1,3})$&$4.3723$&$4.6904$&$K_3\vee4K_1$&$4.6056$&$4.6904$\\
			$K_3\vee(K_1+K_{1,4})$&$6.0322$&$6.4807$&$K_4\vee5K_1$&$6.2170$&$6.4807$\\\hline
		\end{tabular}
}
		\label{tab:2}
\end{table}
From {\it Table 2}, all graphs in $\mathbb{NP}_2$ satisfy $\mu(G)\textless\sqrt{n^2-6n+15}$, a contradiction.

The proof is complete.$\blacksquare$

\noindent{\bf Theorem 3.3} Let $G$ be a connected graph on $n(\ge5)$ vertices with minimum degree $\delta(G)\ge2$. If

\[q(G)\ge\frac{10}{n-1}+2n-6,\]
then $G$ is a pancyclic graph unless $G$ is a bipartite graph or $G=K_3\vee 4K_1$.

\noindent{\bf Proof:} Suppose that $G$ is neither a pancyclic graph nor a bipartite graph. Because $K_n$ is pancyclic and $\delta(K_{1,n-1})=1$. By {\bf Lemma 2.3}
\[\frac{10}{n-1}+2n-6\le q(G)<\frac{2m}{n-1}+n-2,\]
then
\[m>\left(
		\begin{array}{c}
		n-2\\	
		2\\
		\end{array}
	\right)+4.\]

By {\bf Theorem 3.1}, we get $G\in \mathbb{NP}_1$=$\{K_2\vee(K_{n-4}+2K_1), K_5\vee6K_1, K_3\vee(K_2+3K_1), K_3\vee(K_1+K_{1,4}), K_3\vee(K_2+K_{1,3}), (K_2\vee2K_1)\vee5K_1, K_4\vee5K_1, K_{1,2}\vee4K_1, K_2\vee(K_1+K_{1,3}), K_3\vee4K_1\}$. Because when $G\in \{K_2\vee(K_{n-4}+2K_1),  K_3\vee(K_2+K_{1,3}), K_5\vee6K_1\}$, \[m=\left(
		\begin{array}{c}
		n-2\\	
		2\\
		\end{array}
	\right)+4,\] a contradiction.
So, $G\in \mathbb{NP}_2$=$\{K_3\vee(K_1+K_{1,4}), K_3\vee(K_2+K_{1,3}), (K_2\vee2K_1)\vee5K_1, K_4\vee5K_1, K_{1,2}\vee4K_1, K_2\vee(K_1+K_{1,3}), K_3\vee4K_1\}$.

For $(K_2\vee2K_1)\vee5K_1$, let $X=(X_1,X_2,\cdots,X_9)^T$ be the eigenvector corresponding to $q$, where $X_i~(1\le i\le 5)$ correspond to the vertex of degree $4$, $X_i~(6\le i\le 7)$ correspond to the vertex of degree $7$ and $X_i~(8\le i\le 9)$ correspond to the vertex of degree $8$. Then by signless Laplacian eigen-equation (2.2), we have
\begin{equation*}
\setlength{\arraycolsep}{5pt}
\left\{
  \begin{array}{l}
X_1=X_2=\cdots=X_5,X_6=X_7,X_8=X_9,\\
 (q(G)-4)X_1=2X_6+2X_8,\\
 (q(G)-7)X_6=5X_1+2X_8,\\
 (q(G)-8)X_8=5X_1+2X_6+X_8.\\
  \end{array}
\right.
\end{equation*}
Transform above into a matrix equation $(Q'(G)-q(G)I)\tilde X=0$, where $\tilde X=(X_1,X_6,X_8)^T$ and
\begin{equation*}
\setlength{\arraycolsep}{5pt}
Q'(G)=\left[		
		\begin{array}{ccc}		
		4 & 2 & 2\\		
		5 & 7 & 2\\		
		5 & 2 & 9\\
		\end{array}
	\right].
\end{equation*}
Let $g(x) :=det(xI-Q'(G))$, then $g(x)=x^3-20x^2+103x-116$, and $q (G)$ is the largest root of $g(x)=0$. Through calculation, $q(G)=12.5052<\frac{10}{9-1}+2\times9-6=13.2500$ .
Using the same method, we get the signless Laplacian spectral radius of the other graphs in $\mathbb{NP}_2$, showing in the following {\it Table 3},
\begin{table}[H]
\caption{The signless Laplacian spectral radius of $G$}
		\renewcommand\arraystretch{1.5}
		\centering
		\resizebox{\textwidth}{!}{
		\begin{tabular}{C{2.5cm}C{1cm}C{2.5cm}C{2.5cm}C{1cm}C{2.5cm}}\hline
			$G$&$q(G)$&$\frac{10}{n-1}+2n-6$&$G$&$q(G)$&$\frac{10}{n-1}+2n-6$\\\hline
			$K_{1,2}\vee4K_1$&$8.8965$&$9.6667$&$K_3\vee(K_2+K_{1,3})$&$12.6769$&$13.2500$\\
			$K_2\vee(K_1+K_{1,3})$&$9.3408$&$9.6667$&$K_3\vee4K_1$&$9.7720$&$9.6667$\\
			$K_3\vee(K_1+K_{1,4})$&$12.8381$&$13.2500$&$K_4\vee5K_1$&$13.1789$&$13.2500$\\\hline
		\end{tabular}
}
		\label{tab:3}
\end{table}
 From {\it Table 3}, all graphs in $\mathbb{NP}_2$ except $G=K_3\vee 4K_1$ satisfy $q(G)<\frac{10}{n-1}+2n-6$,  a contradiction.

The proof is complete.$\blacksquare$

\section{Appendix}
\begin{figure}[H]
\centering
\includegraphics[width=\textwidth]{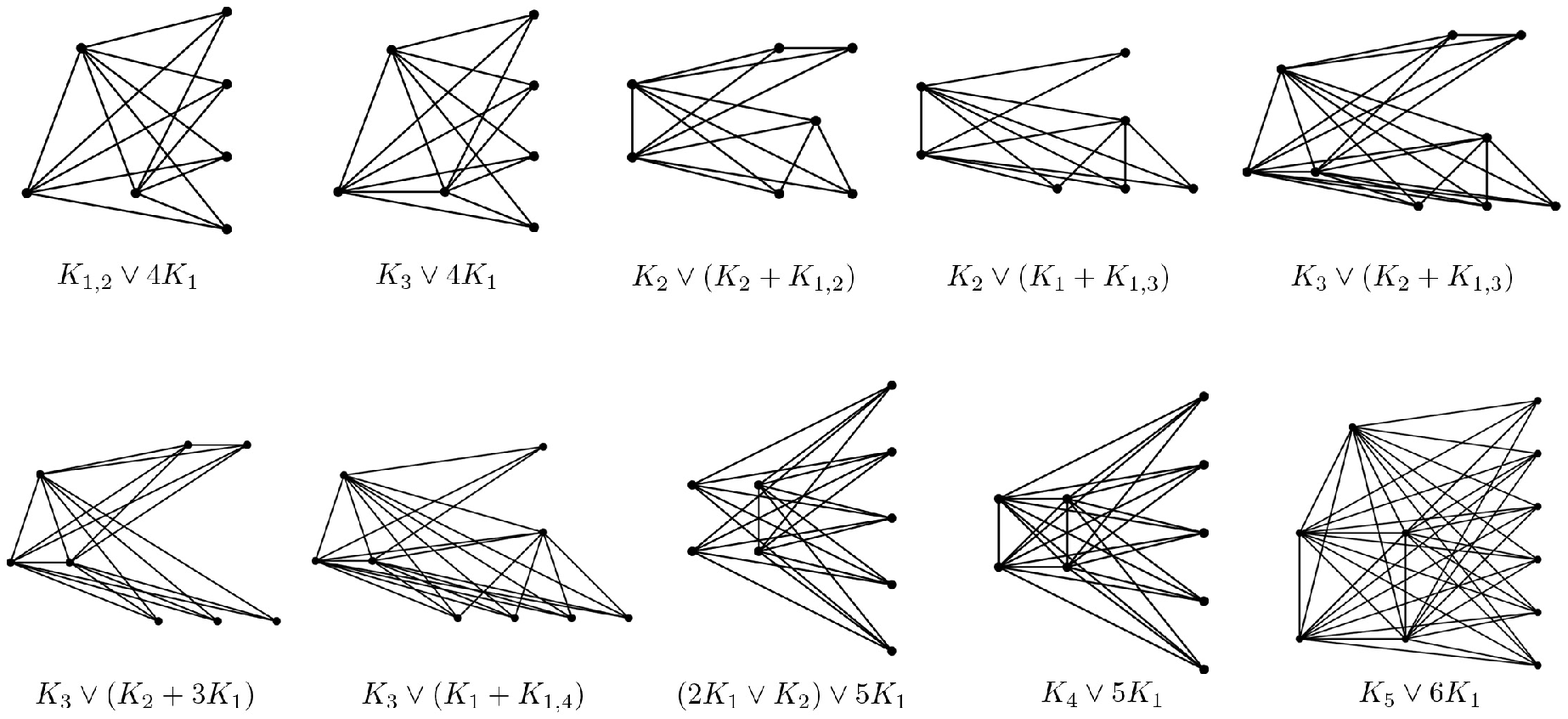}
\includegraphics[width=60mm]{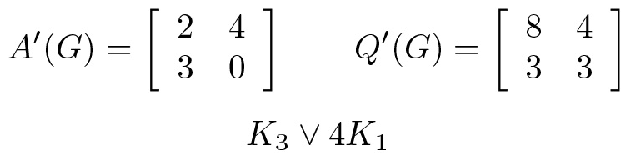}
\includegraphics[width=60mm]{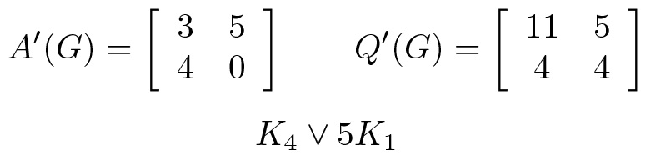}
\includegraphics[width=60mm]{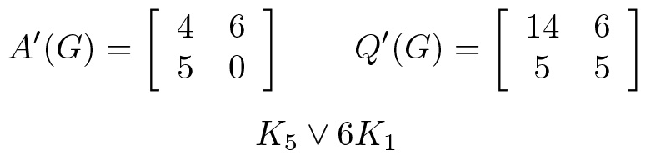}
\includegraphics[width=60mm]{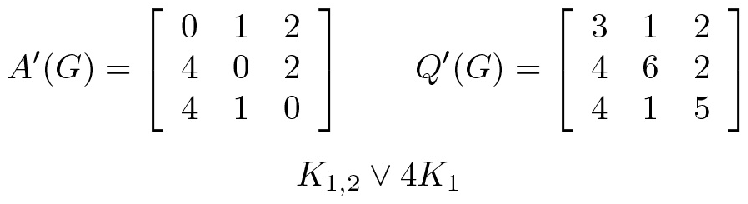}
\includegraphics[width=60mm]{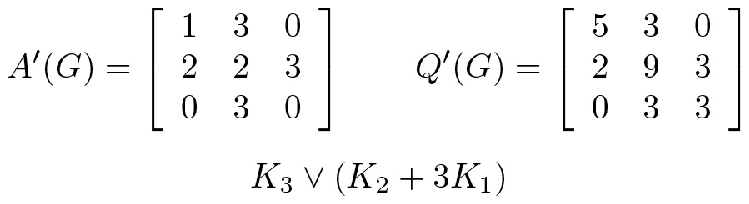}
\includegraphics[width=60mm]{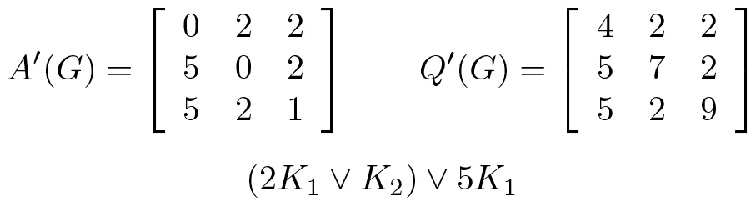}
\includegraphics[width=60mm]{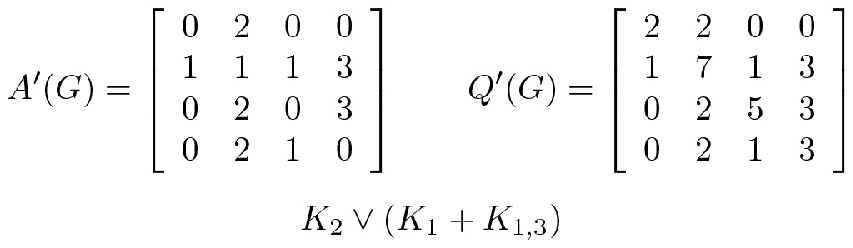}
\includegraphics[width=60mm]{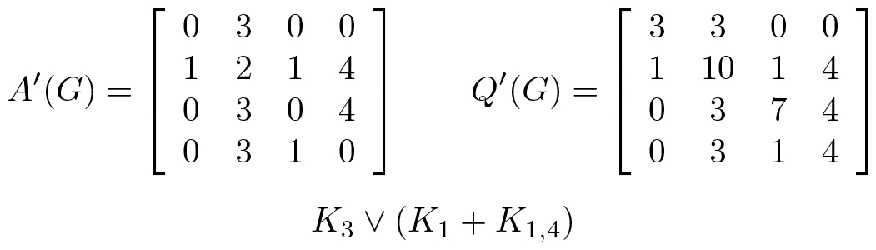}
\includegraphics[width=60mm]{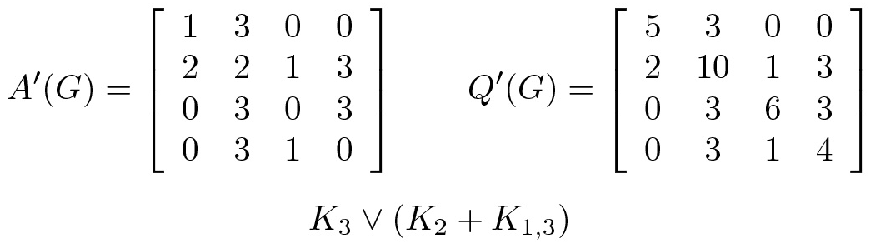}
\includegraphics[width=60mm]{11.eps}
\end{figure}


\end{document}